\documentclass[pmlr]{jmlr} 

\newcommand{\tagdsmode}{proceedings}

\makeatletter
\newcommand{\tagdssubmission}{submission}
\newcommand{\tagdsproceedings}{proceedings}

\ifx\tagdsmode\tagdsproceedings

\else\ifx\tagdsmode\tagdssubmission
  \def\ps@jmlrtps{%
    \let\@mkboth\@gobbletwo
    \def\@oddhead{\scriptsize Under Review at the 2nd Conference on Topology, Algebra, and Geometry in Data Science\hfill}%
    \let\@evenhead\@oddhead
    \def\@oddfoot{}%
    \let\@evenfoot\@oddfoot
  }

\else
  \def\ps@jmlrtps{%
    \let\@mkboth\@gobbletwo
    \def\@oddhead{}%
    \let\@evenhead\@oddhead
    \def\@oddfoot{}%
    \let\@evenfoot\@oddfoot
  }
\fi\fi
\makeatother

\usepackage{longtable}

\usepackage{booktabs}
\usepackage[load-configurations=version-1]{siunitx} 

\theorembodyfont{\upshape}
\theoremheaderfont{\scshape}
\theorempostheader{:}
\theoremsep{\newline}

\jmlrvolume{334}
\jmlryear{2026}
\jmlrworkshop{Topology, Algebra, and Geometry in Data Science}

\title[Scalably computing metric magnitude]{Scalably computing metric magnitude}

\ifx\tagdsmode\tagdssubmission

\else

\author{\Name{Steve Huntsman} \Email{steve.huntsman@cynnovative.com}\\
\AND
\Name{Jewell Thomas} \Email{jewell.thomas@alumni.unc.edu}\\
\AND
\Name{Cynthia Ukawu} \Email{cynthia.ukawu@cynnovative.com}\\
}

\fi

\ifx\tagdsmode\tagdsproceedings
\editor{Editor's name}
\fi

\begin{document}

\maketitle

\begin{abstract}
Applications of metric magnitude often rely on numerically exact results in order to exploit a connection with information theory. We examine various approaches for scaling the dense linear algebra involved and identify hierarchical low-rank solvers as a preferred approach, with a clear path to scales of $10^5$ points on a single powerful workstation, and larger scales using our containerized CUDA-enabled C++/MPI pipeline.
\end{abstract}

\section{\label{sec:introduction}Introduction}

The theory of magnitude (see \S \ref{sec:magnitude} for a quick overview) generalizes notions of size and Euler characteristic to the context of enriched categories \citep{leinster2013magnitude,leinster2016maximizing,leinster2021entropy}. It has already seen remarkable though still preliminary applications in the context of metric spaces, especially of strict negative type, where it generalizes classical information-theoretical constructions to incorporate notions of geometry.

However, applications of metric magnitude have often not scaled, and the techniques have not found wider adoption, primarily because they require solving dense linear systems of the form $Zw = 1$ with $Z_{jk} = \exp(-td_{jk})$ and $d$ a distance matrix. Yet the \emph{weighting vectors} $w$ are exceptionally versatile, with uses for virtually any task involving analysis of Euclidean point clouds and/or information-theoretical characterizations. Accordingly, \citet{andreeva2025approximating} introduced both iterative and greedy schemes for approximating weighting vectors, but neither appears suitable for applications that substantially invoke information-theoretical properties. The iterative scheme produces a nonnegative weighting approximation that is bound to have components with the wrong sign near the most natural and important regime, and the greedy scheme ignores many points without guarantees.

We have therefore undertaken to directly address scalability bottlenecks for numerically exact solutions, in line with work on kernel ridge regression \citep{chavez2020scalable}. We consider three approaches: sparsification, hierarchical low-rank solvers for structured kernel matrices \citep{hackbusch2015hierarchical}, and (detailed in \S \ref{sec:Nystrom}) Nystr\"om approximation. We have determined that while these approaches cannot be fruitfully exploited in conjunction and two of them fail badly, structured kernel matrix solvers show real promise. 

\begin{figure}[htbp]
  \centering
  \includegraphics[trim = 0mm 0mm 0mm 0mm, clip, width=.8\textwidth,keepaspectratio]{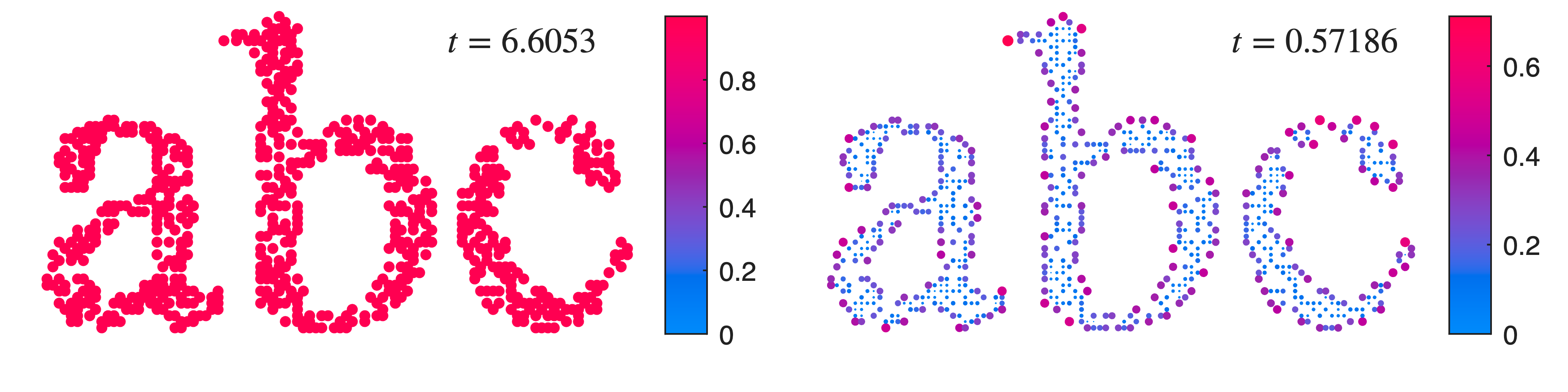}
  \caption{Left: $w$ for $t = \log(n-1)/\min_j \min_{k \ne j} d_{jk} \approx 6.605$, where $d$ is the Euclidean distance on $n = 740$ unique points in $\mathbb{R}^2$ obtained by sampling 1000 points uniformly with replacement from a complicated support over $\mathbb{Z}^2$. Components of $w$ are indicated by size and color: note that all are nearly 1. Right: as in the left panel, with the same points and distance matrix, but for $t = t_+ \approx 0.572$.}
  \label{fig:cutoff}
\end{figure}

There is an easy bound $t_+ \le \log(n-1)/\min_j \min_{k \ne j} d_{jk}$ for the interval $[t_+,\infty)$ where a natural information-theoretical interpretation of $w$ is available \citep{huntsman2023diversity}, but determining $t_+$ precisely requires solving $Zw = 1$ for different values of $t$. In this paper, we focus on this problem (and retaining $w$ as a byproduct) for data in $\mathbb{R}^m$. Figure \ref{fig:cutoff} shows $w$ for both the bound above and for $t_+$ itself on the same underlying points.
\footnote{
For $t \in (0,t_+)$, the diversity-maximizing distribution can be $\mathbf{NP}$-hard to compute, since it requires consideration of exponentially many sparsity patterns. Moreover, while in Euclidean space there is a polynomial-time algorithm for $t = 0$ \citep{huntsman2025peeling}, it still requires computing $d^{-1}1$ for a series of distance submatrices, which is at least as hard as computing $Z^{-1}1$ for $t > 0$.
}

The upside of scalably solving $Zw = 1$ is highlighted by demonstrations and applications to, e.g., classification, active learning, and anomaly detection \citep{bunch2020practical}; computer vision \citep{adamer2024magnitude}; sparse reinforcement learning \citep{guo2021geometric}; introspection of neural representations \citep{bunch2021weighting,andreeva2023metric,limbeck2024metric}; pooling in graph neural networks \citep{limbeck2026geometry}; dataset quality/curation \citep{couch2024beyond} with implications for network distillation; and single- and multi-objective optimization \citep{huntsman2022parallel,huntsman2023diversity,huntsman2023quality,pereverdieva2025comparative,emmerich2026magnitude}, and even chemistry \citep{bi2024persistent,bi2025topological}. In short, many aspects of machine learning could benefit from techniques built on weighting vectors that were computed more scalably. 

Our contributions are i) a systematic empirical analysis of scaling approaches for solving $Zw = 1$; ii) a containerized, open-source C++/MPI pipeline built on STRUMPACK \citep{strumpack2020} 
at \url{https://github.com/Cynnovative/compass},
and iii) documented failures of other hitherto plausible scaling approaches that may be useful to other researchers.

The paper is organized as follows. \S \ref{sec:magnitude} gives background on magnitude, weightings, and diversity; \S \ref{sec:scalingIntro} discusses approaches to scaling; \S \ref{sec:methods} details our methods and experiments; and \S \ref{sec:scaling} discusses scaling. Appendices \S \ref{sec:sparsificationExamples} and \S \ref{sec:Nystrom} respectively give examples illustrating failures of sparsification and Nystr\"om approximations; and \S \ref{sec:semianalytical} details various semianalytical results that may be of interest but fail to produce obvious algorithmic advantage.

\section{\label{sec:magnitude}Weightings, magnitude, and diversity}

We now get into details. A square matrix $Z \ge 0$ with $\text{diag}(Z) > 0$ is called a \emph{similarity matrix}. We consider similarity matrices of the form $Z = \exp[-td]$ where $(f[M])_{jk} := f(M_{jk})$ indicates componentwise function application, $t \in (0,\infty)$, and $d$ is a square matrix with entries in $[0,\infty]$ that also satisfy the triangle inequality. In this paper, we only consider Euclidean metric matrices of the form $d_{jk} = \|x_{(j)} - x_{(k)}\|$ for $\{x_{(j)}\}_{j=1}^n \subset \mathbb{R}^m$.


A \emph{weighting} $w$ is a solution to
\begin{equation}
\label{eq:weighting}
Zw = 1, 
\end{equation}
where $1$ indicates a vector of all ones. If $Z$ has a weighting $w$, its \emph{magnitude} is $\text{Mag}(Z) := 1^T w$. If $d$ is Euclidean, then $Z$ is positive definite and it has a unique weighting. Weightings turn out to be excellent scale-dependent boundary or outlier detectors in Euclidean space \citep{willerton2009heuristic,bunch2020practical} due to a link with Bessel capacities \citep{meckes2015magnitude}. 

\begin{example}
\label{ex:3PointSpace}
Let $\{x_j\}_{j=1}^3 \subset \mathbb{R}^2$ with $d_{jk} := d(x_j,x_k)$ given by $d_{12} = d_{13} = 1 = d_{21} = d_{31}$ and $d_{23} = \delta = d_{32}$ with $\delta \ll 1$. A brief calculation shows that
\[w_1 = \frac{e^{(\delta+2)t}-2e^{(\delta+1)t}+e^{2t}}{e^{(\delta+2)t}-2e^{\delta t}+e^{2t}}; \quad w_2 = w_3 = \frac{e^{(\delta+2)t}-e^{(\delta+1)t}}{e^{(\delta+2)t}-2e^{\delta t}+e^{2t}}.\]
As Figure \ref{fig:3pointSpace} shows, for $t \ll 1$, $w \approx (1/4,1/4,1/2)^T$; for $t \gg 1$, $w \approx (1,1,1)^T$, and for $t \approx 10$, $w \approx (1/2,1/2,1)^T$. Thus the ``effective number of points'' goes from $\approx 1/4+1/4+1/2 = 1$, to $\approx 1/2 + 1/2 + 1 = 2$, to $\approx 1+1+1 = 3$. 


\begin{figure}[h]
  \centering
  \includegraphics[trim = 10mm 0mm 10mm 0mm, clip, width=.8\textwidth,keepaspectratio]{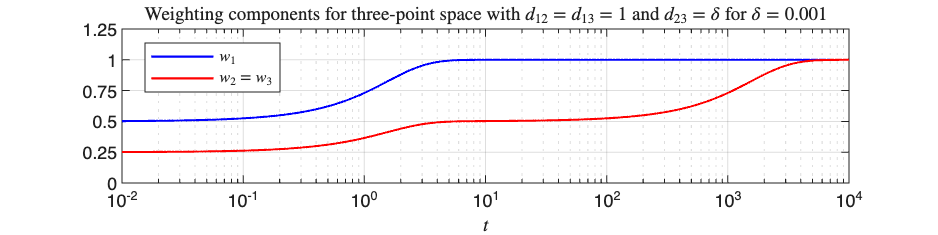}
  \caption{The magnitude $w_1+w_2+w_3$ is a scale-dependent ``effective number of points.''}
  \label{fig:3pointSpace}
\end{figure}
\end{example}

Magnitude and weightings play a key role \emph{maximizing} a general and axiomatically defined measure of diversity \citep{leinster2016maximizing,leinster2021entropy}. 
For $1 < q < \infty$, the \emph{diversity of order $q$} for a probability distribution $p$ and similarity matrix $Z$ is defined as
\begin{equation}
\label{eq:diversity}
\log D_q^Z(p) := \frac{1}{1-q} \log \sum_{j: p_j > 0} p_j (Zp)_j^{q-1}
\end{equation} 
and extended via limits for $q = 1,\infty$. $D_q^Z(p)$ is a ``correct'' measure of diversity in very much the same way that Shannon entropy is a ``correct'' measure of information. In fact, $\log D_q^Z(p)$ is a geometrical generalization of the R\'enyi entropy of order $q$. The usual R\'enyi entropy corresponds to $Z = I$, and Shannon entropy corresponds to $Z = I$ and $q = 1$.

\begin{theorem}
\label{thm:maxDiversity}
If $Z$ is symmetric, positive definite, and has a unique positive weighting $w$, then for all $q$, $w$ is proportional to $\arg \max_p D_q^Z(p)$ \citep{leinster2016maximizing}.
\end{theorem}

Theorem \ref{thm:maxDiversity} reduces diversity maximization to linear algebra while rendering $q$ irrelevant. Accordingly, the practically canonical nonzero scale for computing weightings is the \emph{cutoff} $t_+$, i.e., the least value of $t$ such that $w \ge 0$: see Figure \ref{fig:cutoff}. The cutoff can be computed efficiently via root-finding for $\min_j w_j$ using the elementary bound $t_+ \le \log(n-1)/\min_j \min_{k \ne j} d_{jk}$. In practice Ridders' method \citep{ridders1979new,press2007numerical} performs better than bisection due to its prioritization of low function evaluation count rather than low iteration count.
\footnote{It is useful to halt the root-finder when $\min_j w_j \in (0, \varepsilon)$ for $\varepsilon$ comparable to machine epsilon. This has the additional advantage of producing nondegenerate maximum-diversity probability distributions.}

The cutoff $t_+$ is the least scale that permits the application of Theorem \ref{thm:maxDiversity}. However, in practice the cutoff scale is often quite large in the sense that the resulting weighting has many components with values close to unity, limiting its utility. It is often more advantageous to work in the limit $t = 0$, though the corresponding weightings are often sparsely supported. This limit is examined in \citet{huntsman2025peeling} and \citet{huntsman2026peel}. Here, we focus on finding $t_+$, though structured/kernel solvers can also be used for the $t = 0$ case.

\section{\label{sec:scalingIntro}Approaches to scaling}

There are several promising approaches to solve $Zw = 1$ more efficiently at scale: e.g., enforcing sparsity, or using linear solvers for structured kernel matrices. 
\S \ref{sec:Nystrom} and \S \ref{sec:semianalytical} respectively detail failures of Nystr\"om approximations and semianalytical techniques.

\subsection{\label{sec:sparsificationIntro}Sparsification}

The obvious approach to sparsity is to threshold the similarity matrix $Z$, implicitly setting small entries to zero, say by efficient approximate $K$-nearest neighbor search using hierarchical navigable small world graphs \citep{malkov2018efficient}. However, the theory in its present state crucially requires $Z$ to be positive definite to have useful guarantees about any information-theoretical interpretation of $w$. While this is guaranteed for a Euclidean distance matrix by classical results, thresholding $Z$ usually destroys positive definiteness \citep{guillot2015functions,guillot2016preserving}. While taking $t$ large enough ensures positive definiteness even after thresholding, approximation errors still result, and in practice we want $t$ as close to $t_+$ as possible, which in turn requires solving $Zw = 1$ accurately. 

\begin{figure}[htbp]
  \centering
  \includegraphics[trim = 10mm 1mm 10mm 1mm, clip, width=.8\textwidth,keepaspectratio]{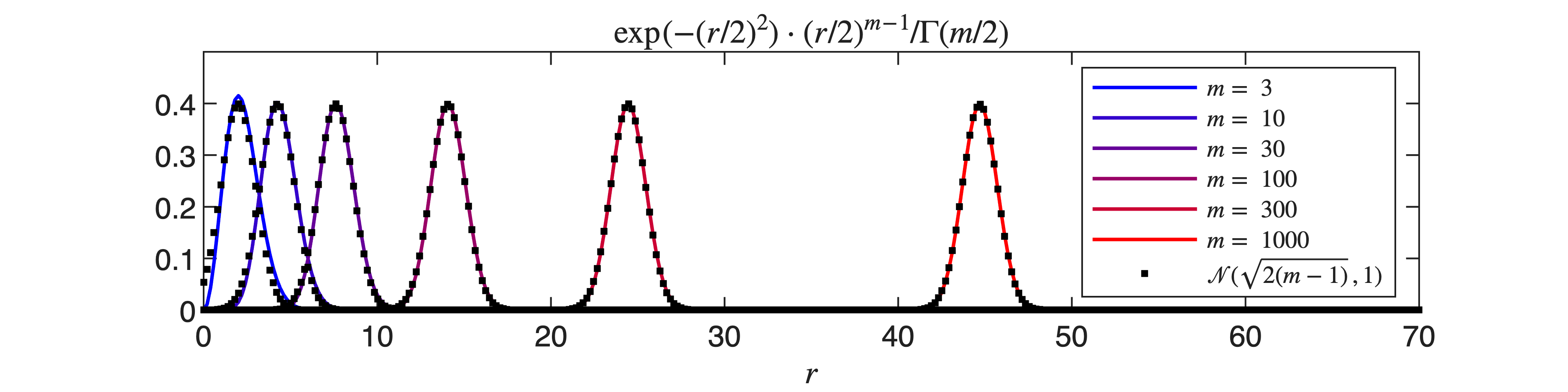}
  \caption{The distribution 
  of distances between two IID points $\sim \mathcal{N}(0,I_m)$.}
  \label{fig:gaussianDistanceDistribution}
\end{figure}

Meanwhile, ``distance concentration'' often effectively makes different pairwise distances indiscernible in high dimension. The probability distribution of distance $d_{12} = r$ between two IID points $\sim \mathcal{N}(0,I_m)$ 
rapidly approximates $\mathcal{N}(\sqrt{2(m-1)},1)$, as Figure \ref{fig:gaussianDistanceDistribution} shows.
\footnote{
To derive the distribution in Figure \ref{fig:gaussianDistanceDistribution},
note that the difference of two IID $\mathcal{N}(0,I_m)$ points is distributed as $\mathcal{N}(0,2 I_m)$: the norm of this difference is chi-distributed (in fact, $\chi_m(r/\sqrt{2})$). A mechanical derivation is in \citep{thirey2015distribution}, and a more sophisticated partial generalization to $\mathcal{N}(0,\Sigma)$ with $\Sigma_{jj} \equiv 1$ is in \citep{martirosyan2024euclidean}. Meanwhile, a formula for the distribution of distances between two independent uniformly sampled points from a convex body is in \citep{aharonyan2020distribution}.
}
While the essentially constant width of the distribution produces a narrowing of normalized distances, Gaussian tail behavior also raises the prospect that a quasilinear number of point pairs might control $Z$. However, numerical and analytical examinations of Gaussian-distributed data show that sparsification is not a tenable route to fast and accurate solutions of $Zw = 1$ at scale. Even as a first step, any effective sparsification via thresholding would require the identification of a small number of so-called ``hubs'' with many reverse nearest neighbors \citep{angiulli2018behavior}, which is also challenging. But in any event our numerical results from sparsification were so unsatisfactory as to not warrant significant discussion: see \S \ref{sec:sparsificationExamples}.

    

\subsection{\label{sec:structuredIntro}Solvers for structured kernel matrices}

If $d$ is a Euclidean distance matrix, then the corresponding $Z$ is automatically positive definite, so we get a reproducing kernel Hilbert space to which Mercer’s theorem applies. In practice this amounts to the existence of a Cholesky decomposition. But a more promising method to solve $Zw = 1$ at scale is to use efficient hierarchical low-rank solvers that exploit kernel structure in a more detailed way and offer subquadratic performance. However, like the fast multipole method that spawned them, hierarchical solvers traditionally scale very badly with ambient problem dimension, or at best with data dimension \citep{march2015kernel}.

Our experiments show that these techniques can offer real advantage: while dense kernel solvers scale well with dimension,  hierarchical solvers permit much larger scale applications, albeit with dimension dependence.

\section{Methods, experiments, and results}\label{sec:methods}

Our numerical experiments were performed on a 16-core MacBook Pro M4 Max with 128 GB RAM and ARM CPU (the integrated GPU was not exercised since it did not offer CUDA). We implemented two main computational pipelines: a high-performance C++/MPI backend for large kernel matrix experiments using STRUMPACK \citep{strumpack2020}, and a Python-based pipeline for dense and sparse approaches. All experiments are conducted in Docker containers to ensure reproducibility of environment configurations. Docker images are configured to run on either x86 or ARM CPUs. We made all 16 CPUs available for Docker and capped its memory at 66.5 GB. Thread and device environments (\texttt{OMP\_NUM\_THREADS}, \texttt{OPENBLAS\_NUM\_THREADS} etc.) are set in the appropriate containers. 


\subsection{Software environment}

\paragraph{C++/STRUMPACK pipeline.}
For large-scale experiments involving hierarchical low-rank approximation and kernel solves, we developed a C++ driver compiled against the STRUMPACK library (v7.2.0), with MPI parallelism enabled and GPU support. The current build uses OpenBLAS (v0.3.21) and OpenMPI for BLAS/LAPACK and parallel support, with optional CUDA (12.3+) for GPU acceleration. SLATE (v2025.05, for distributed BLAS) and METIS (v5.1) were also linked for fast matrix reordering and factorization.

\paragraph{Python/SciPy pipeline.}
For dense solving, we used \texttt{scipy\hspace{0pt}.linalg\hspace{0pt}.solve} in Python 3.12, with multithreaded BLAS operations supported by OpenBLAS. All package versions are managed using \texttt{micromamba} and are specified via YAML file. 

\paragraph{Docker configurations.}
Our C++ pipeline is based on \texttt{nvidia\hspace{0pt}/cuda:12.3.1\hspace{0pt}-devel\hspace{0pt}-ubuntu22.04} to enable leveraging CUDA GPU in the future. Our Python pipeline is based on \texttt{mambaorg\hspace{0pt}/micromamba:2.3.1} to ensure we can link SciPy and NumPy to the most performant linear algebra libraries (e.g., OpenBLAS) and leverage multi-threaded processing.

\subsection{Experimental workflow}


For both the C++ and Python pipelines, the key algorithm is root-finding for $t_+$, as described above. As a sanity check that adds marginal time, we test positive (semi)definiteness by estimating the smallest eigenvalue (C++: symmetric Lanczos via STRUMPACK; Python: \texttt{scipy.sparse.linalg.eigsh}). To speed up the Python pipeline, we first perform a fast Cholesky check with \texttt{scipy.linalg.cholesky} and only perform an eigenvalue decomposition if this check fails. We solve \eqref{eq:weighting} using (\texttt{strumpack::HSSMatrix::solve}) in C++ and \texttt{scipy.linalg.solve} in Python.


All code, build recipes, environment files, and Docker images (for both C++ and Python stacks) are in GitLab and will be publicly released to ensure full reproducibility of all results. 

\subsection{Experimental results}

To demonstrate the pipeline and gauge runtimes, we sampled $n$ points from $\mathcal{N}(0,I_m)$ and from $\mathcal{N}(0,\Sigma)$ with $\Sigma_{jk} := \delta_{jk} 2^{-(j-1)}$. We take $n \in \{1000,3000,10000,30000\}$ and $m \in \{3,10,30,100\}$, and we perform three runs per configuration.
\footnote{In general, the computational complexity of solving is at least linear in $m$ because of the underlying kernel evaluations. If the intrinsic dimension of the data is fixed while the ambient dimension $m$ grows, hierarchical solvers will generally continue to perform well \citep{march2015kernel}.}
For each set of points, the computational task is to compute $t_+$ and the weighting at that scale, as discussed above. 

In preliminary experiments not detailed here, sparsification resulted in unacceptable performance, with Kendall tau coefficients relative to exact results overwhelmingly below $0.85$. Also, on relatively small problems, iterative methods (either through SciPy or STRUMPACK) solved as fast or faster than the dense pipeline, with similar levels of fidelity, but scaling behavior suggested that iterative methods are only marginally beneficial.

\begin{figure}[htbp]
  \centering
  \includegraphics[trim = 0mm 0mm 0mm 2mm, clip, width=\textwidth,keepaspectratio]{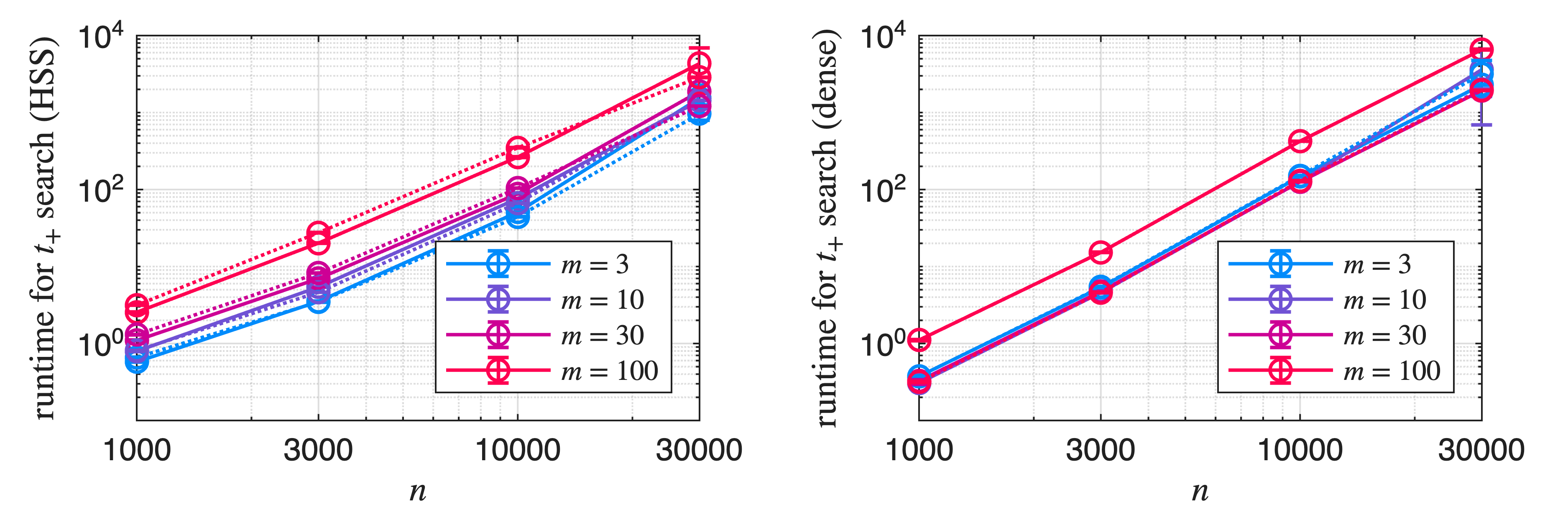}
  \caption{Left: total runtime (in seconds) for $t_+$ search using HSS methods. Results for $\Sigma = I_m$ are shown as solid lines; results for $\Sigma_{jk} := \delta_{jk} 2^{-(j-1)}$ are shown as dashed lines. Error bars are obtained from three experiments for each pair $(m,n)$. Right: as in the left panel, but for a dense solver. Not shown: the number of search iterations for $m = 100$ was much higher for the dense solver.}
  \label{fig:time_s}
\end{figure}

\begin{figure}[htbp]
  \centering
  \includegraphics[trim = 0mm 0mm 0mm 0mm, clip, width=\textwidth,keepaspectratio]{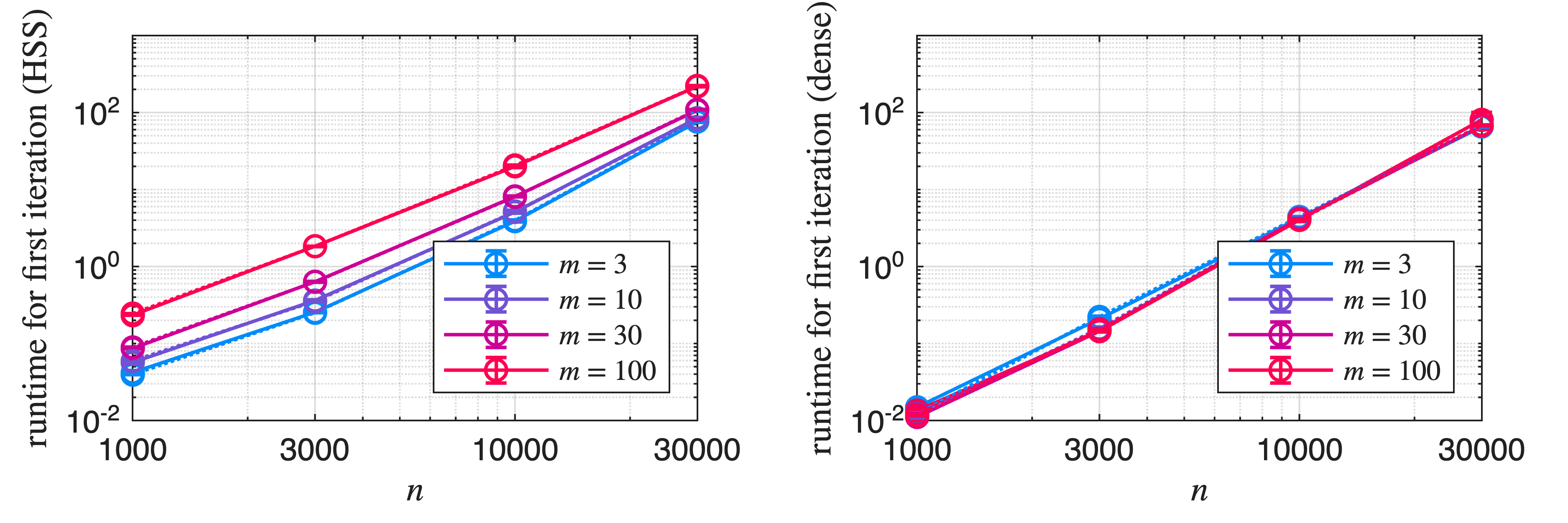}
  \caption{As in Figure \ref{fig:time_s}, but for runtime (in seconds) for the first iteration of $t_+$ search.}
  \label{fig:first_iter_s}
\end{figure}

\begin{figure}[htbp]
  \centering
  \includegraphics[trim = 0mm 0mm 0mm 2mm, clip, width=\textwidth,keepaspectratio]{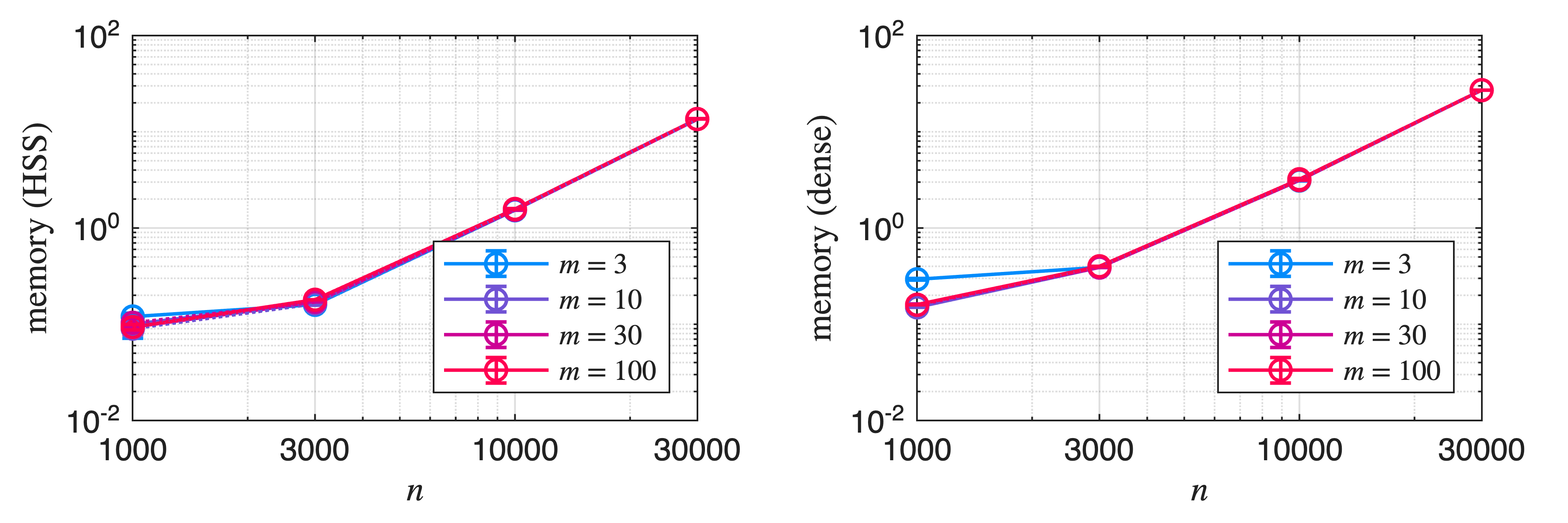}
  \caption{As in Figure \ref{fig:time_s}, but for memory usage (in GB).}
  \label{fig:mem_gb}
\end{figure}

\begin{figure}[htbp]
  \centering
  \includegraphics[trim = 0mm 0mm 0mm 2mm, clip, width=\textwidth,keepaspectratio]{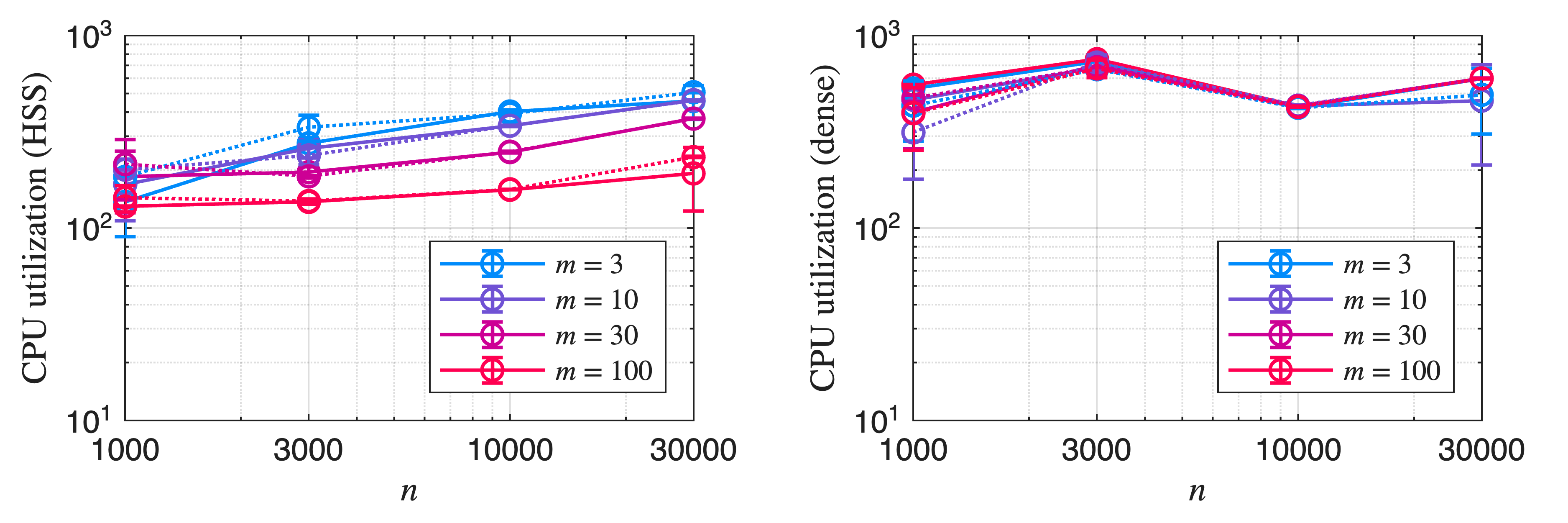}
  \caption{As in Figure \ref{fig:time_s}, but for CPU utilization percentage.}
  \label{fig:cpu_pct}
\end{figure}

Figures \ref{fig:time_s}-\ref{fig:cpu_pct} show that HSS methods through STRUMPACK scale better than ordinary dense methods from SciPy (even with BLAS optimization and multithreading): both runtime and memory use are noticeably better for the former for $n \gtrsim 10^4$, with significantly lower CPU utilization above low dimension. Because STRUMPACK is built for distributed computation, it is reasonable to visually extrapolate our time and memory results, provided that runtime is appropriately amortized over nodes. Notably, although HSS performance degrades with dimension (as expected), this degradation is more graceful than a catastrophic impact that might have been expected.

\section{Scaling}\label{sec:scaling}

\subsection{Current scaling}\label{sec:currentScaling}

Our experiments deliberately focused on laptop-scale performance: however, more scaling is clearly possible using distributed MPI on a multi-node cluster and/or CUDA \citep{wang2019exact}. We have implemented this and tested it, albeit on a single machine without CUDA, where overhead 
hinders performance. However, this is the correct approach for large problems on a multi-node cluster where MPI ranks have access to additional cores and memory bandwidth. Our Dockerfile builds with CUDA support, so access to NVIDIA hardware can provide GPU acceleration as well, primarily due to an internal $ULV$ factorization involving dense block operations on the HSS tree.

For $n \approx 30000$, the time for a single iteration of the rootfinding algorithm (including but not just solving a single instance of $Zw = 1$) with STRUMPACK is typically on the order of one minute, which is comparable to MATLAB.
\footnote{
There is no way to reuse factorizations of $Z(t)$ for different values of $t$, which informs the problem of finding the cutoff $t_+$.
}
Matrix compression takes about a quarter of that time; Lanczos slightly under ten percent; and $ULV$ factorization takes about two thirds of the time. Around $n \approx 70000$, our solver segfaults inside OpenBLAS's \texttt{dgemm\_itcopy} during compression, as a debugger backtrace confirms. This is an issue with BLAS, not STRUMPACK, and the practical ceiling for our current build is $n = 65000$. The underlying cause is likely large HSS ranks and/or compression tolerances.

We attempted various optimizations with little or no effect. For example, reducing the number of Lanczos iterations from 100 to 30 provided real but marginal performance gains. Loose HSS tolerances led to unacceptable noise, and we avoided them. Relaxing compression tolerances does not meaningfully reduce HSS compression ranks, because Gaussians are hard to compress. Replacing scattered heap allocations (\texttt{std::vector<std::vector<double>>}) with a single array improved cache locality and let the compiler vectorize distance computation with SIMD instructions, which accelerated compression by several percent, but compression itself is only about a quarter of the overall runtime.

\subsection{Future scaling}\label{sec:futureScaling}

Any frontal assault on large-scale instances of \eqref{eq:weighting} requires STRUMPACK or a functional analogue (e.g., GOFMM \citep{yu2017geometry}) built for supercomputers. While scales of $n \approx 10^5$ are achievable on powerful workstations, substantially larger scales require clusters
at $n \approx 10^6$, or leading supercomputing resources at $n \approx 10^8$. 

Unfortunately, faster exact solutions or good approximations are not clearly in evidence, but many practical applications may not require either. For example, a black-box global optimizer that uses weighting vectors under the hood \citep{huntsman2022parallel,hoffman2022benchmarking} might not actually require precise solutions of \eqref{eq:weighting} in order to efficiently explore a function landscape. Good local approximations and/or bad global approximations would be very scalable and might work well enough even for hundreds of thousands or millions of points in hundreds or thousands of dimensions. We leave this, and particularly the investigation of alternative low-rank approximations, for future work.

\acks{
Thanks to Evan Gorman for many useful conversations; and to Yang Liu for providing advice regarding STRUMPACK; and to referees for suggestions that helped the presentation.

We used Claude Opus to help find references and evaluate implementations.

This research was partially developed with funding from the Defense Advanced Research Projects Agency (DARPA). The views, opinions and/or findings expressed are those of the authors and should not be interpreted as representing the official views or policies of the Department of Defense or the U.S. Government. 
}

\bibliography{scaling}

\appendix

\section{\label{sec:sparsificationExamples}Examples of sparsification for Gaussian-distributed points}

\begin{figure}[htbp]
  \centering
  \includegraphics[trim = 20mm 30mm 20mm 25mm, clip, width=\columnwidth,keepaspectratio]{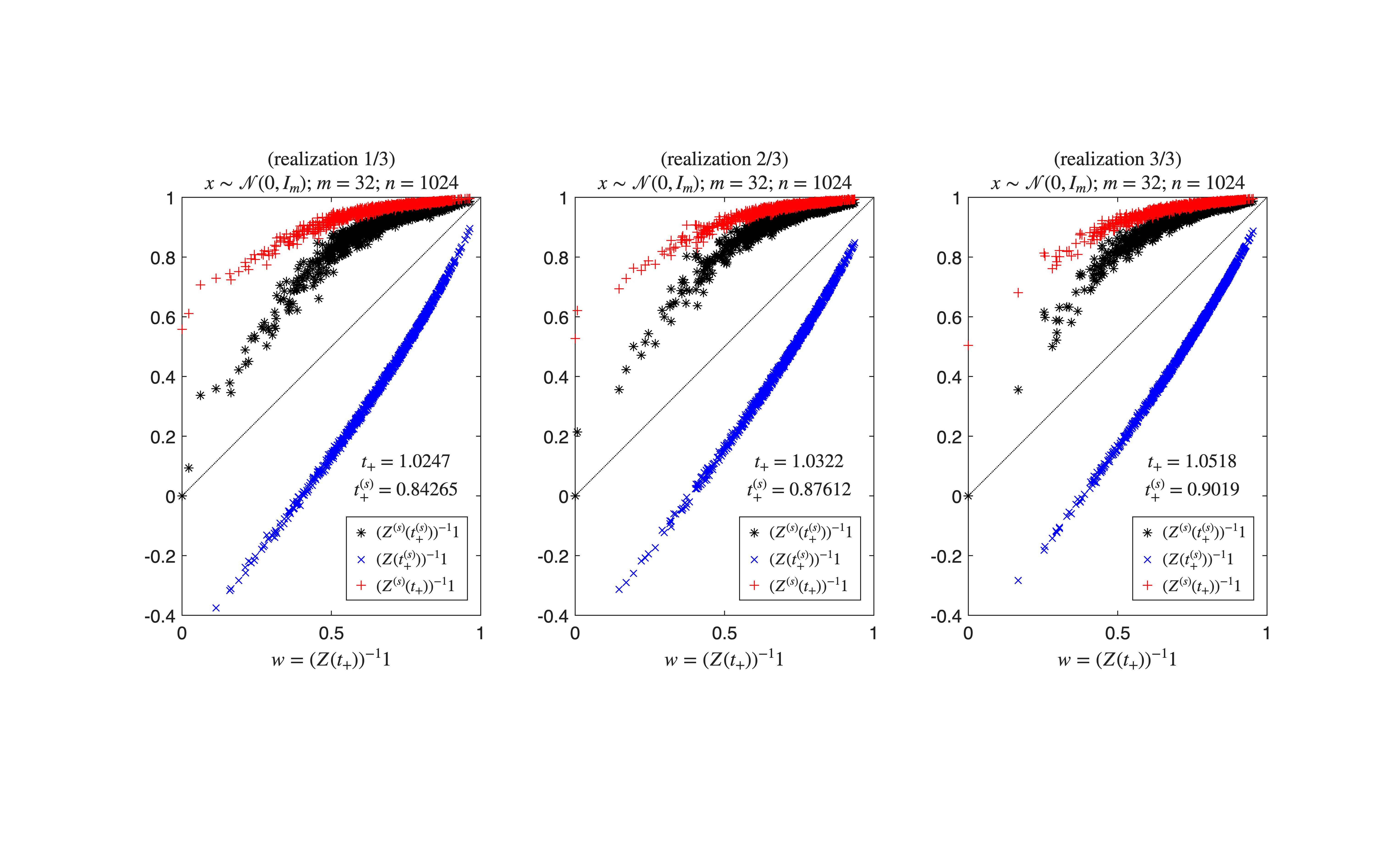}
  \caption{Effects of sparsification for three realizations of a sample of size $n = 1024$ from $\mathcal{N}(0,I_m)$ with dimension $m = 32$. For each sample, we compute the true cutoff $t_+$ and corresponding similarity matrix $Z$, as well as a symmetrized sparse approximation $Z^{(s)}$ using $\lceil \log_2(n) \rceil = 10$ approximate nearest neighbors and approximate cutoff $t_+^{(s)}$ computed using the same binary search technique as for $t_+$. We show weighting approximations using {\color{blue}the sparse cutoff (blue $\times$ markers)}, {\color{red}the sparse similarity matrix (red $+$ markers)}, and both simultaneously (black $*$ markers).}
  \label{fig:gaussian_sparsification}
\end{figure}

Figure \ref{fig:gaussian_sparsification} illustrates how sparsification using $K$ approximate nearest neighbors obtained as in \citep{malkov2018efficient} (and still symmetrizing) severely degrades solutions to $Zw = 1$ and to the calculation of cutoffs for Gaussian-distributed points. The bad approximations shown are broadly typical, though much worse results occurred frequently in other realizations, while materially better results did not occur.

It may be surprising that sparsification has numerical effects as profound as Figure \ref{fig:gaussian_sparsification} indicates. However, it is unsurprising in light of \S \ref{sec:sparsificationIntro} that sparsification destroys many practical theoretical guarantees as well.

\section{\label{sec:Nystrom}Nystr\"om approximation}

The Nystr\"om approximation \citep{cai2022fast} for a kernel matrix is of the form
\begin{equation}\label{eq:Nystrom}
\begin{pmatrix}
    K_{(LL)} & K_{(LM)} \\ K_{(LM)}^T & K_{(MM)}
\end{pmatrix} 
\approx
\begin{pmatrix}
    K_{(LL)} & K_{(LM)} \\ K_{(LM)}^T & K_{(LM)}^T K_{(LL)}^+ K_{(LM)}
\end{pmatrix}
\end{equation}
where a partition of points into a small ``landmark'' set $L$ and a large ``main'' set $M$ is indicated. It is reasonable to hope that this might help scalably approximate $Z^{-1}1$ and $d^{-1}1$. 
\footnote{
Although Nystr\"om approximations are useful for preconditioning \citep{frangella2023randomized,abedsoltan2024nystrom}, we do not consider this in detail here since STRUMPACK has robust preconditioning facilities.
}
However, naive applications are essentially unusable for our (or similar) purposes that demand precision in order to compute $t_+$ and the associated weighting.

The most immediate reason is that the kernels $\exp(-t\|x-x'\|)$ and $\|x-x'\|$ generically produce matrices without any good low-rank approximation (cf. Figures \ref{fig:random_kernel_matrix_m128n512} and \ref{fig:random_kernel_matrix_m128n512_scale0}). This leads to a numerical catastrophe. For the exponential kernel, this is easily mitigated by subtracting the diagonal of the right hand side of \eqref{eq:Nystrom} and adding an identity matrix. However, this mitigation is not just obviously imperfect, but it also leaves some room for improvement.

Specifically, the resulting approximate weighting vectors tend to behave particularly badly at landmark points. Using uniformly random landmarks and performing affine transformations on the landmark and main components to match their first two moments while holding the main variance fixed yields an approximate weighting vector that is highly correlated with the exact result for a wide range of scales. However, the approximation is only quantitatively good at one scale, and it is not even clear how to find this scale \emph{a priori}. 

Using a greedy stochastic algorithm to find more ``outlying'' landmarks along the lines of Algorithm 2 in \citep{huntsman2023quality}
requires a different approach for transforming landmark and main components. In this event, it is more appropriate to transform the main components as above, but to transform the landmarks differently, by matching them with the upper tail of the main components, then shifting the landmark components by the average difference between this matching transformation and the one above. 

\begin{figure}[htbp]
  \centering
  \includegraphics[trim = 45mm 5mm 45mm 5mm, clip, width=.8\textwidth,keepaspectratio]{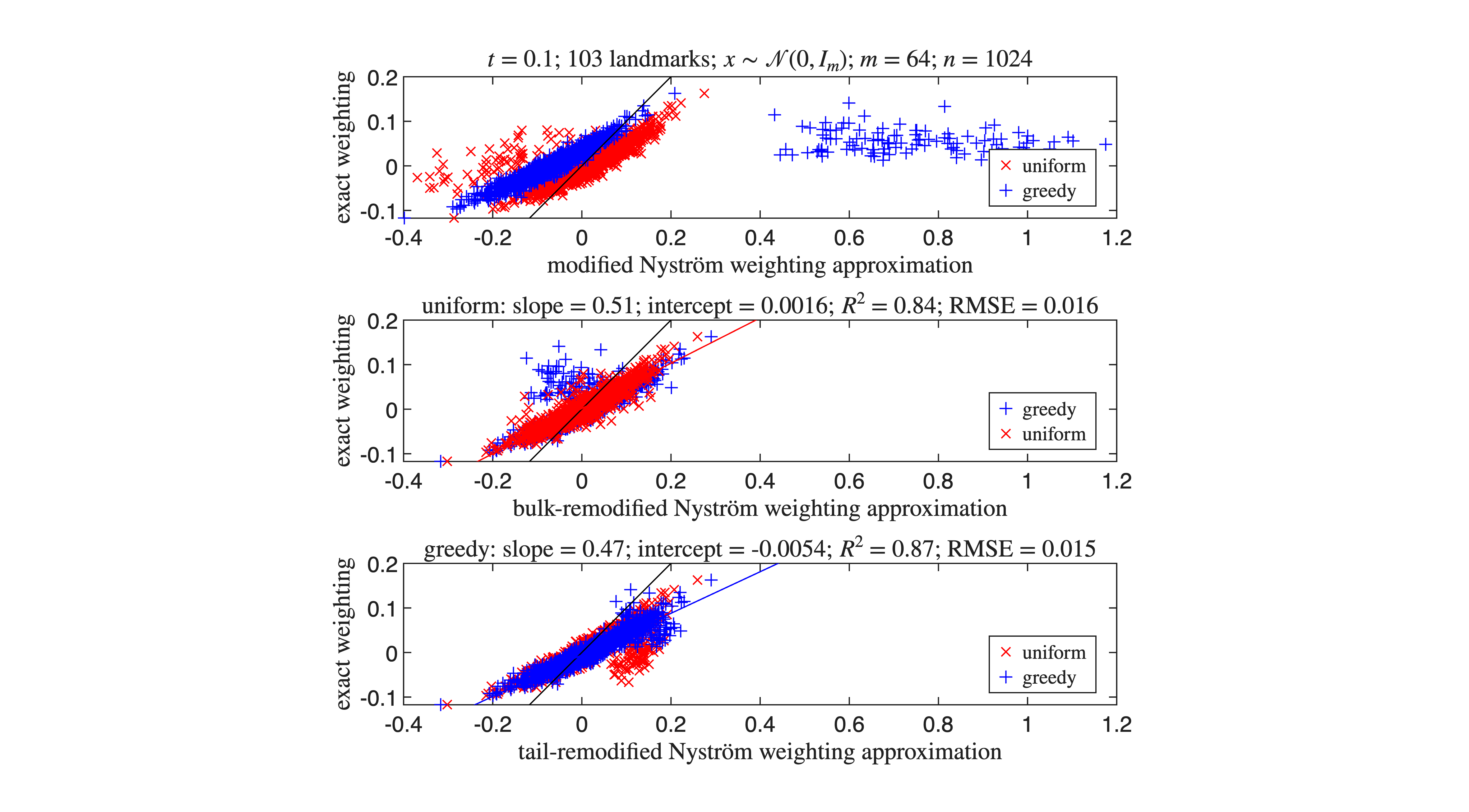}
  \caption{A comparison of modified Nystr\"om approximations to a weighting using {\color{red}uniform} and {\color{blue}stochastic magnitude-greedy} landmarks with $t = 0.1$. The upper panel compares the approximations obtained with a unit diagonal to the exact result, shown as the black diagonal line. In each case, the outliers are the landmarks. The middle panel shows the results of affine transformations on landmark and main components by basic moment-matching: this works better for uniform landmarks: a corresponding linear fit is also shown. The lower panel shows the results of a more delicate transformation suitable for greedy landmarks: a corresponding linear fit is also shown.}
  \label{fig:nystrom_10}
\end{figure}

\begin{figure}[htbp]
  \centering
  \includegraphics[trim = 50mm 5mm 50mm 5mm, clip, width=.32\textwidth,keepaspectratio]{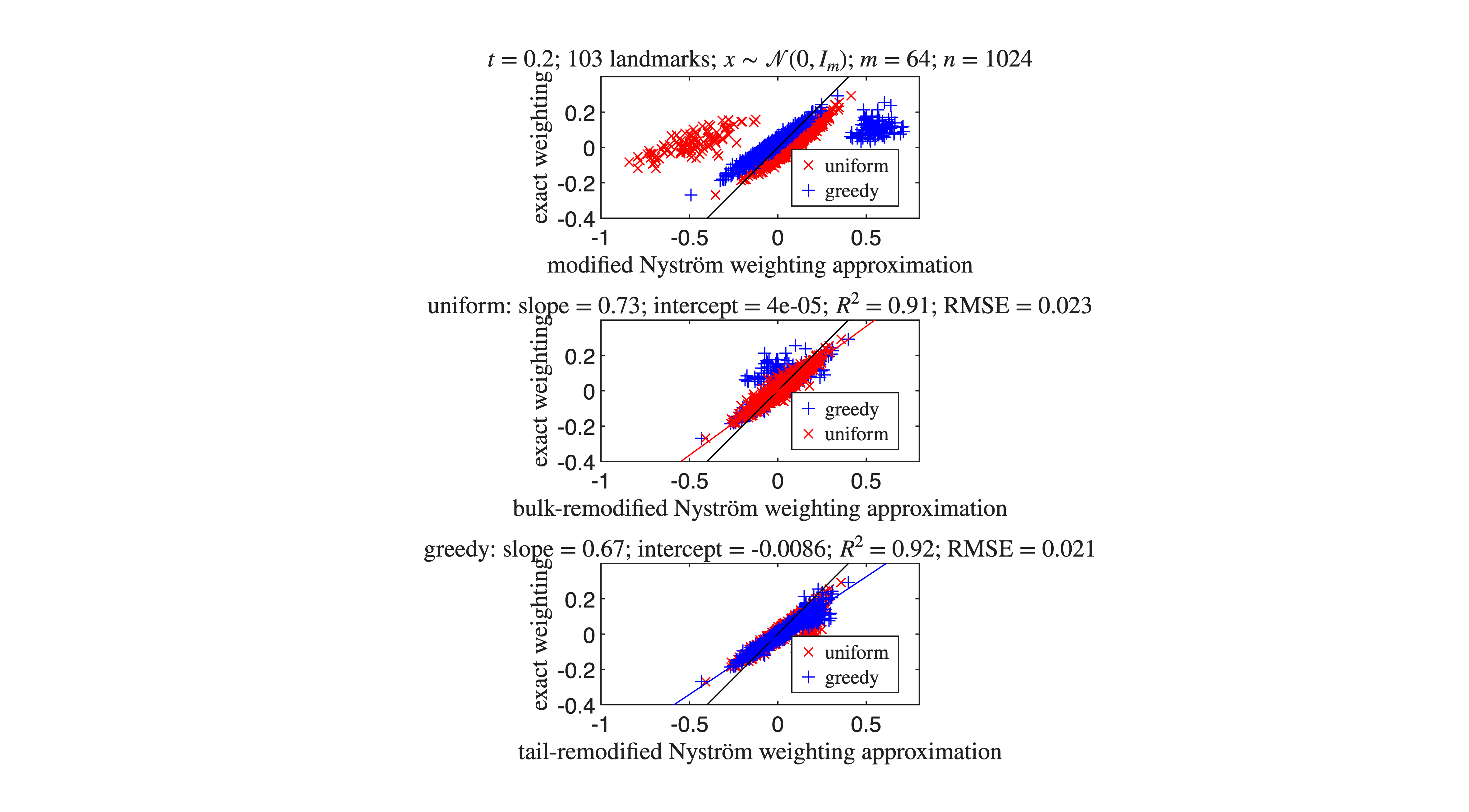}
  \includegraphics[trim = 50mm 5mm 50mm 5mm, clip, width=.32\textwidth,keepaspectratio]{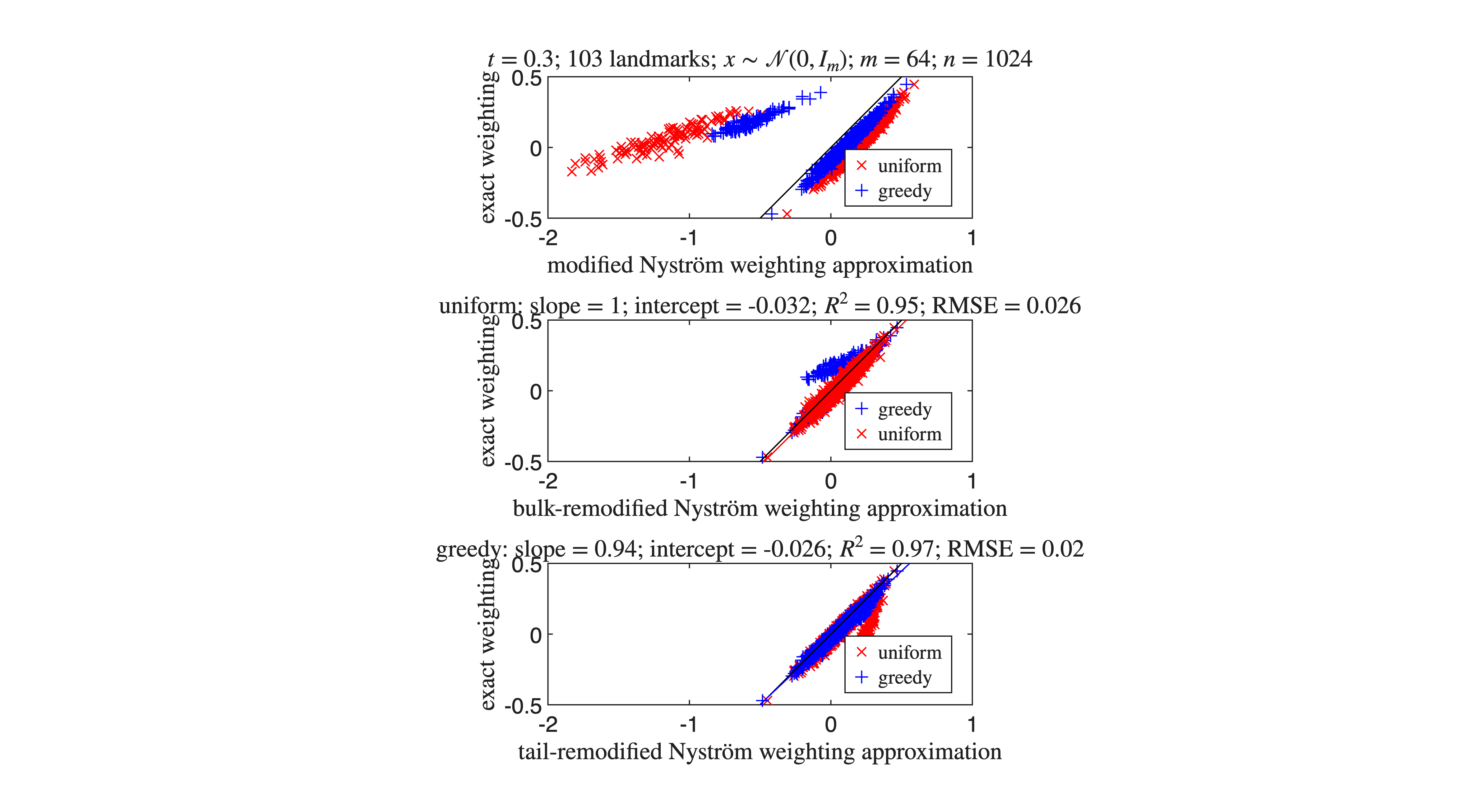}
  \includegraphics[trim = 50mm 5mm 50mm 5mm, clip, width=.32\textwidth,keepaspectratio]{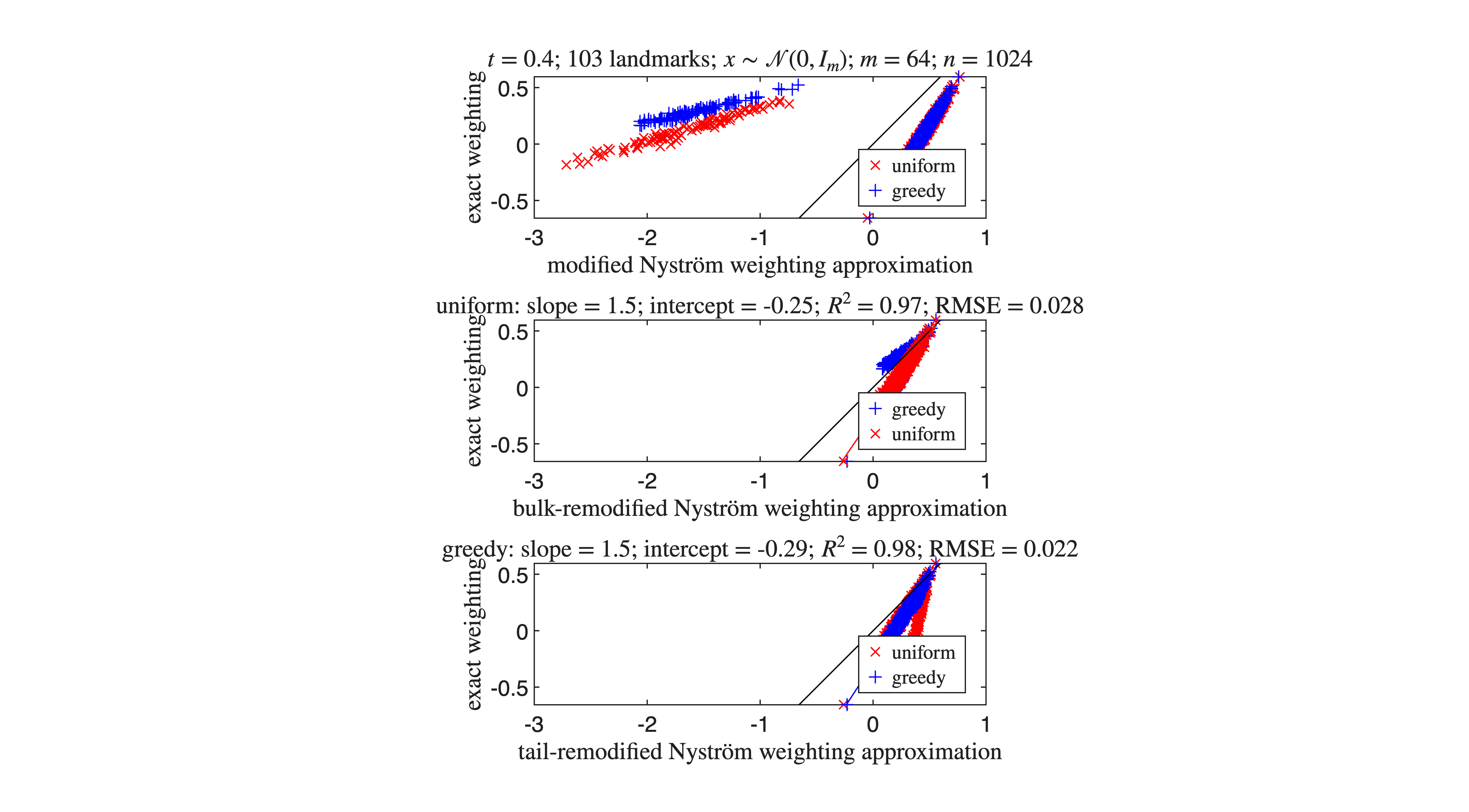}
  \caption{As in Figure \ref{fig:nystrom_10}, but for $t = 0.2$ (left), $t = 0.3$ (center), and $t = 0.4$ (right).}
  \label{fig:nystrom_20and30and40}
\end{figure}

This more careful application of Nystr\"om approximations holds more promise than a naive application, but it still has no obvious utility. Figures \ref{fig:nystrom_10} and \ref{fig:nystrom_20and30and40} show how these approximations are strongly correlated with exact results, but with uncertain linear fits that strongly depend on $t$. Unfortunately, even if we can use this approximation to reliably identify the index of the minimal weighting component, it is not clear how to use that information to compute the cutoff scale more efficiently. For example, the intercepts of linear fits when $t = 0.4$ are far from zero, while the actual cutoff scale is just under $0.7$ (cf. Figures \ref{fig:gaussian_cutoff_part_two} and \ref{fig:gaussian_cutoff_part_one}).

\section{\label{sec:semianalytical}Semianalytical investigations}

\subsection{\label{sec:semianalyticalWeightingScaleZero}Weightings at scale zero for Gaussian-distributed points}

In practice, weightings at scale zero are often more useful than weightings at finite nonzero scale \citep{huntsman2025peeling}. While the latter requires solving $Zw = 1$ once for $Z = \exp[-td]$ and $d$ a Euclidean (or more generally here, strict negative type) distance matrix, the former requires solving equations of the form $dw = 1$, since
\begin{equation}\label{eq:scaleZero}
    \lim_{t \downarrow 0} (\exp[-td])^{-1}1 = \frac{d^{-1}1}{1^T d^{-1} 1}.
\end{equation}
To get diversity-maximizing distributions, it is necessary to solve equations of the form $dw = 1$ repeatedly, albeit at progressively smaller scales.

Meanwhile, a ``canonically hard'' case for solving the weighting equation at zero or finite scale is presented by a sample of $n$ IID points $x^{(j)}$ from a standard $m$-dimensional Gaussian, i.e., $x^{(j)} \sim \mathcal{N}(0,I_m)$ for $j \in [n]$. 

For the scale zero case and recalling the approximation $\mathcal{N}(\sqrt{2(m-1)},1)$ of the distribution of distances between IID points $\sim \mathcal{N}(0,I_m)$, we have that the linear equation $dw = 1$ is approximated by the equation 
\begin{equation}\label{eq:Wigner1}
(A+\sqrt{2(m-1)}(11^*-I))w = 1,
\end{equation} 
where $A_{jk} \sim \mathcal{N}(0,1)$ for $j \ne k$ and $A_{jj} \equiv 0$. 
\footnote{\label{foot:finiteScaleWigner}
The case for finite scale is obviously better behaved, but harder to work with analytically and also less useful: $Z_{jk} \sim \text{LogNormal}(-t\sqrt{2(m-1)}(1-\delta_{jk}),t^2)$. 
}
An IID instance of \eqref{eq:Wigner1} is a so-called Wigner ensemble \citep{akemann2011oxford}, but to our knowledge the structure of solutions to $dw=1$ has not been studied. Our case of course has entries that are obviously not independent even if approximately identically distributed. 

As we shall see, the structure of solutions to an instance of \eqref{eq:Wigner1} with IID upper triangular entries of $A$ differs dramatically from \eqref{eq:weighting} for Gaussian-distributed points. But let us entertain the alternative hope briefly. As justification, it is reasonable to try to characterize and compare solutions to the Wigner ensemble and actual scale zero weighting equations (for which entries are not independent) with the same entrywise statistics to try to develop a heuristic for sparsification. Because the standard Gaussian mass is concentrated in a thin spherical shell, the distance matrix behaves particularly badly from the point of view of sparsification, and any heuristic that worked effectively in this situation would be likely to be applicable very generally. 
\footnote{\label{foot:frame}
A stringent test of such a heuristic would be an overcomplete frame that is approximately equiangular and tight, say by perturbing overcomplete equiangular tight frames such as those enumerated in \citep{fickus2015tables}. This would arguably be the worst possible case to confront in practice, though the prospects for characterization are too weak to warrant starting here.
}
    
Applying the Sherman-Morrison formula to \eqref{eq:Wigner1} generically yields
\begin{equation}\label{eq:WignerSM}
\left ( A+\sqrt{2(m-1)}(11^*-I) \right )^{-1}1 = \frac{\left ( A - \sqrt{2(m-1)}I \right )^{-1}1}{1+1^* \left ( A - \sqrt{2(m-1)}I \right )^{-1}1} 
\end{equation} 
and
\begin{align}\label{eq:WignerPerturbation}
\left ( A - \sqrt{2(m-1)}I \right )^{-1}1 & = -\frac{1}{\sqrt{2(m-1)}} \left ( I - \frac{1}{\sqrt{2(m-1)}} A \right )^{-1}1 \nonumber \\
& \approx -\frac{1}{\sqrt{2(m-1)}} \left ( I + \frac{1}{\sqrt{2(m-1)}} A \right )1.
\end{align}
In the IID case, $A1 \sim \mathcal{N}(0,(n-1)I)$, so the maximum-diversity distribution is approximately uniform plus a Gaussian perturbation. 
\footnote{\label{foot:highOrderWigner}
The next higher order correction to \eqref{eq:WignerPerturbation} is governed under the IID \emph{Ansatz} by
\begin{align}
(A^2 1)_j & = \sum_{k\ell} A_{jk} A_{k\ell} \nonumber \\
& = \sum_{k,\ell \ne j} A_{jk} A_{k\ell} +  \sum_k A_{jk}^2 \nonumber \\
& \sim \left ( \pi^{-1} K_0(| \cdot |) \right )^{*^{(n-1)^2}} * \chi^2_{n-1} \nonumber
\end{align}
and in principle (perhaps even practice) this can be evaluated further. Higher-order corrections are yet more complicated. Even the expectations of $A^\ell 1$ under the IID \emph{Ansatz} involve Isserlis' theorem \citep{munthe2025short} as a building block. Meanwhile, the actual products of IID Gaussian random variables are distributed as Meijer $G$-functions \citep{springer1966distribution,stojanac2017products} whose tail behavior has been analyzed in \citep{leipus2023distribution}.
}
The approximation \eqref{eq:WignerPerturbation} could conceivably not be as bad as it looks: while the entries of the left hand side are more often negative than not, the denominator of the factor multiplying it is distributed as $\mathcal{N}(1-n/\sqrt{2(m-1)},n(n-1)/4m^2)$, and $1-n/\sqrt{2(m-1)}$ will typically be negative in practice. 

\begin{figure}[htbp]
  \centering
  \includegraphics[trim = 0mm 0mm 0mm 0mm, clip, width=\columnwidth,keepaspectratio]{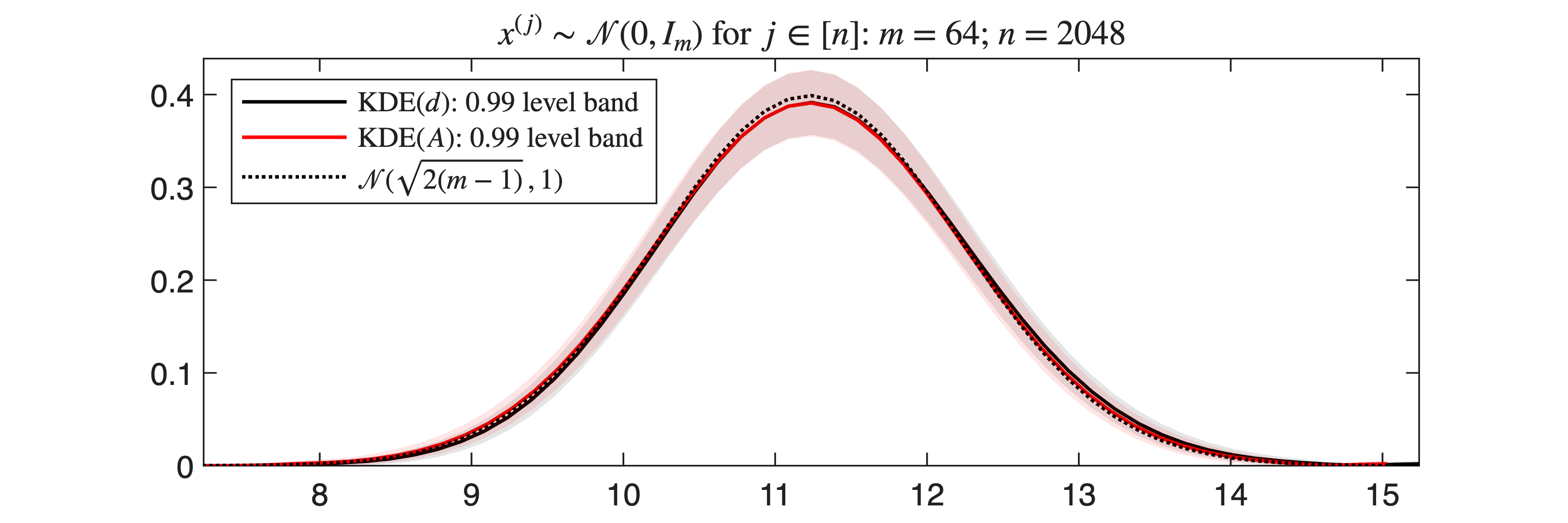}
  \caption{Density estimates for $d$ and $A$ along the lines of \eqref{eq:Wigner1}.}
  \label{fig:semianalytical_attempt_part_one}
\end{figure}

\begin{figure}[htbp]
  \centering
  \includegraphics[trim = 0mm 0mm 0mm 0mm, clip, width=\columnwidth,keepaspectratio]{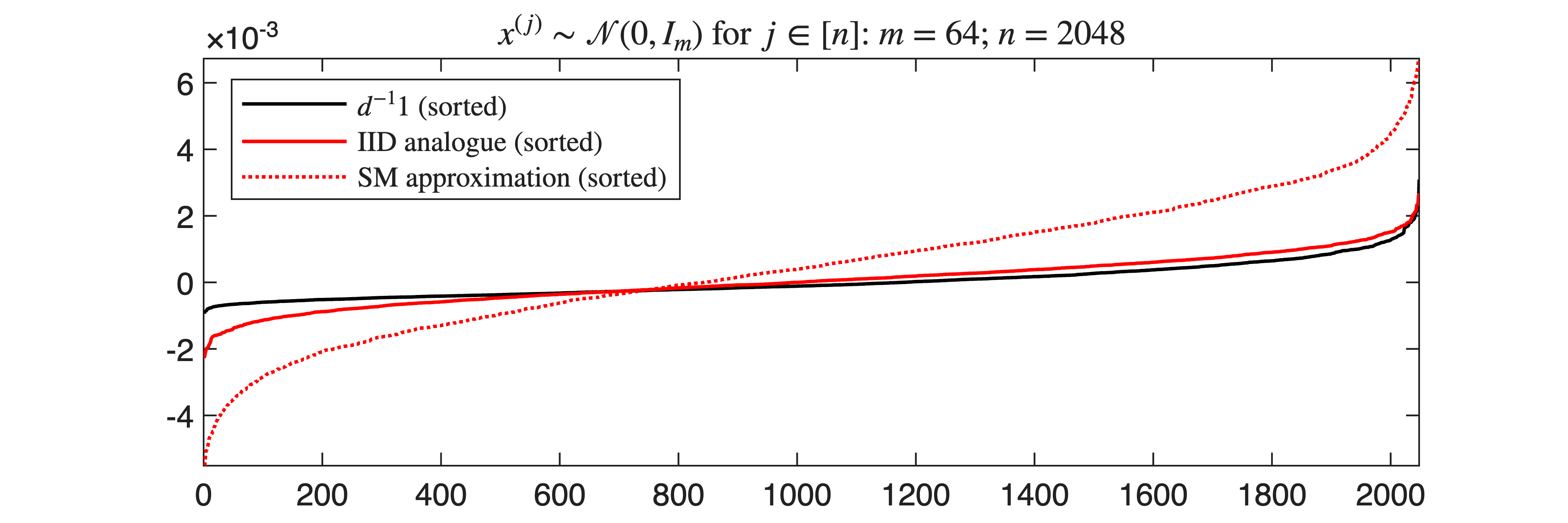}
  \caption{Sorted entries of $d^{-1}1$, the IID analogue \eqref{eq:Wigner1}, and a subsequent approximation to \eqref{eq:WignerSM} using \eqref{eq:WignerPerturbation}.}
  \label{fig:semianalytical_attempt_part_two}
\end{figure}

In any event, numerics are required to get a practical understanding and to see if this model has any use at all. And it turns out that numerics quickly dispel the hope entertained above. Figure \ref{fig:semianalytical_attempt_part_one} illustrates (kernel density estimates of) the entries of $d$ and $A$ along the lines of \eqref{eq:Wigner1}, while Figure \ref{fig:semianalytical_attempt_part_two} illustrates how the IID approximation \eqref{eq:Wigner1} to $d^{-1}1$ and its further approximation \eqref{eq:WignerSM} are very bad.

\subsection{\label{sec:semianalyticalCutoff}Cutoff scales for Gaussian-distributed points}

\begin{figure}[htbp]
  \centering
  \includegraphics[trim = 0mm 0mm 0mm 0mm, clip, width=\columnwidth,keepaspectratio]{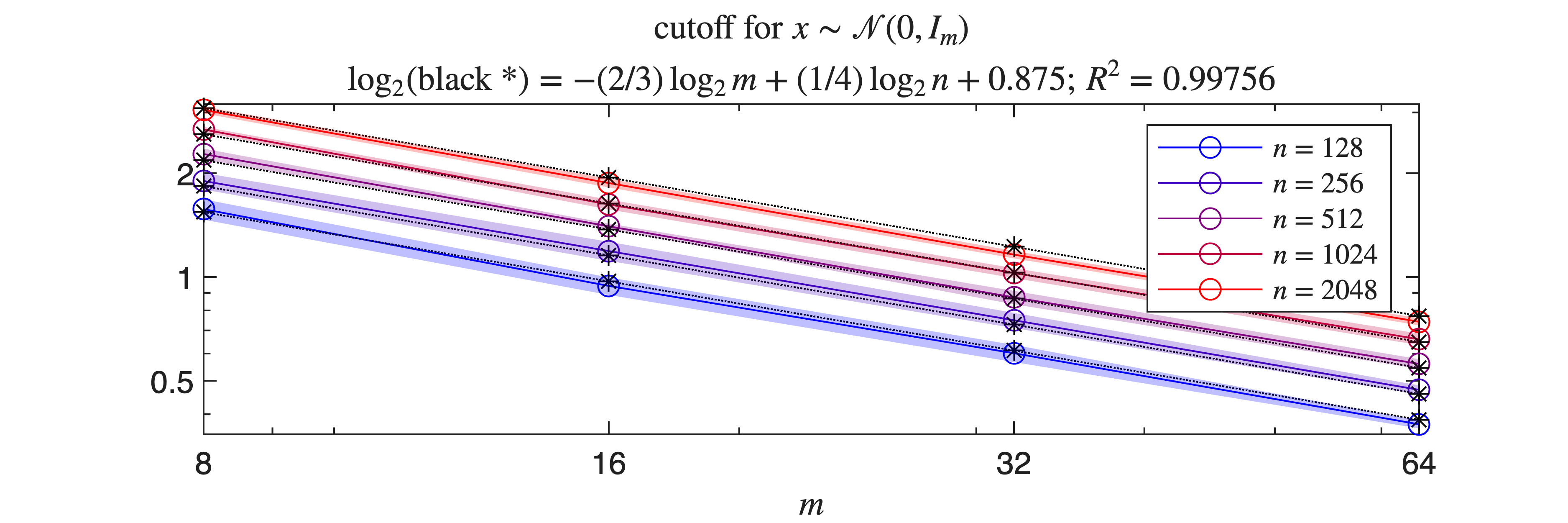}
  \caption{Cutoffs of $N = 10$ realizations of size $n$ from $\mathcal{N}(0,I_m)$, with standard deviations indicated, and with \eqref{eq:gaussianCutoff} in black for comparison (taking $\gamma_+ \approx 2^{7/8} \approx 1.834$. For smaller $m$ and $n$ the linear behavior breaks down and the standard deviations are much larger.}
  \label{fig:gaussian_cutoff_part_two}
\end{figure}

\begin{figure}[htbp]
  \centering
  \includegraphics[trim = 0mm 0mm 0mm 0mm, clip, width=\columnwidth,keepaspectratio]{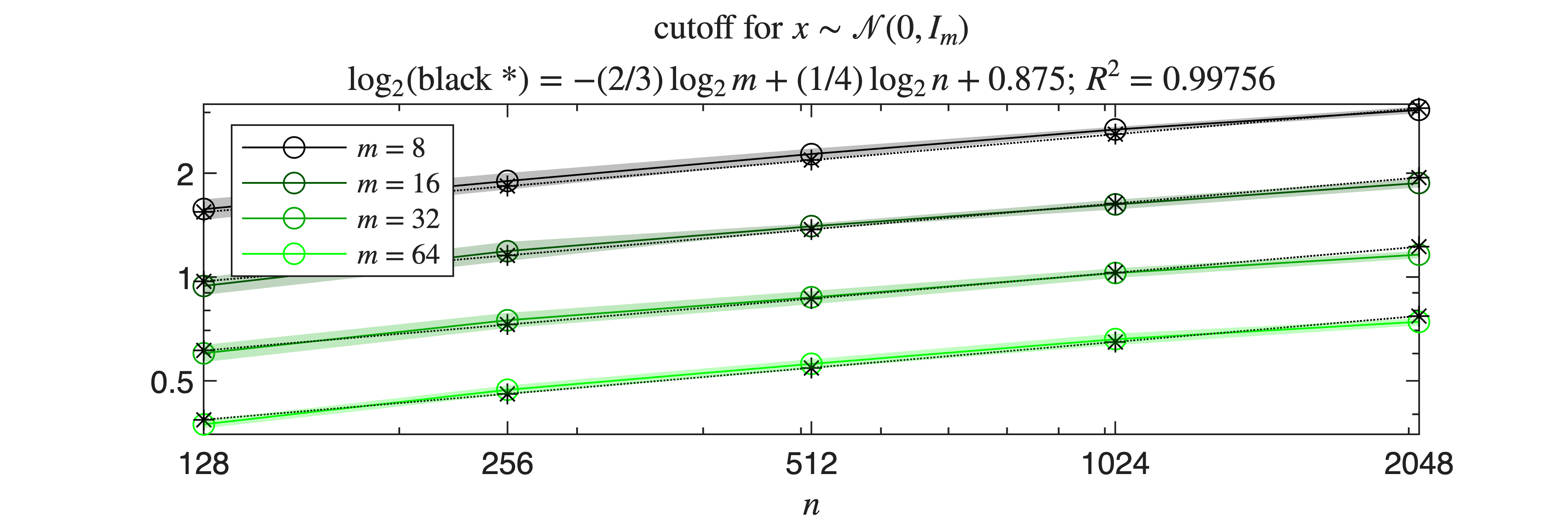}
  \caption{As in Figure \ref{fig:gaussian_cutoff_part_two}, but as a function of sample size $n$ vs. of dimension $m$.}
  \label{fig:gaussian_cutoff_part_one}
\end{figure}

Intriguingly, a numerical experiment suggests that for $x^{(j)} \sim \mathcal{N}(0,I_m)$ with $j \in [n]$, we have the excellent approximation
\begin{equation}\label{eq:gaussianCutoff}
\mathbb{E}t_+ \approx \gamma_+ m^{-2/3} n^{1/4}
\end{equation}
where $\gamma_+$ is slightly less than 2. Figures \ref{fig:gaussian_cutoff_part_one} and \ref{fig:gaussian_cutoff_part_two} show the results of computing cutoffs for $n$ points distributed as $\mathcal{N}(0,I_m)$ for varying dimension $m$ and sample size $n$.

It is natural to try to exploit \eqref{eq:gaussianCutoff} to compute actual cutoffs more quickly. One idea along these lines is to use \eqref{eq:gaussianCutoff} to produce an initial guess for a cutoff that might be better than the easy upper bound. The trivial scaling relationship $t_+(\lambda d) = \lambda^{-1} t_+(d)$ indicates how to normalize $d$ towards this end. However, this idea does not seem to save (or cost) much time, and it is not clear how to exploit \eqref{eq:gaussianCutoff} in practice. 

Still, \eqref{eq:gaussianCutoff} is sufficiently intriguing to warrant further investigation.

\subsection{\label{sec:semianalyticalRandom}Random kernel matrices}

Let $x^{(j)} \sim \mathcal{N}(0,m^{-1}I_m)$ be IID for $j \in [n]$, and as usual let $d_{jk} = \|x^{(j)} - x^{(k)}\|$ be the corresponding distance matrix. This is a special case of a \emph{Euclidean random matrix} \citep{karoui2010spectrum,bordenave2013euclidean}, or more generally, a \emph{random kernel matrix} \citep{do2013spectrum} of the form 
\begin{equation}\label{eq:randomKernelMatrix}
M_{jk} := F(x^{(j)},x^{(k)},m). 
\end{equation}
Random kernel matrices concentrate around their mean \citep{amini2021concentration}, and as such they tend to have bad numerical behavior: similarity matrices are no exception here.

On the other hand, Theorem 2 of \citep{do2013spectrum} gives the nice result that for $$F(x^{(j)},x^{(k)},m) = f(\|x^{(j)} - x^{(k)}\|^2)$$ with $f$ differentiable at $2$, $M$ has the same limiting spectral distribution as 
\begin{equation}\label{eq:sameSpectralDistributionAsRKM}
B := \left [ f(0) - f(2) + 2f'(2)\right ] I_n - 2f'(2) X^TX,
\end{equation}
where $X_{j\ell} := x^{(j)}_\ell$. Since $\text{spec}(\alpha I + \beta A) = \alpha + \beta \cdot\text{spec}(A)$, the spectral distribution of \eqref{eq:sameSpectralDistributionAsRKM} follows trivially from the spectral distribution for the Wishart ensemble $X^T X$, which tends to a Marchenko-Pastur distribution \citep{bai2010spectral,akemann2011oxford}.

\begin{figure}[htbp]
  \centering
  \includegraphics[trim = 0mm 0mm 0mm 0mm, clip, width=\textwidth,keepaspectratio]{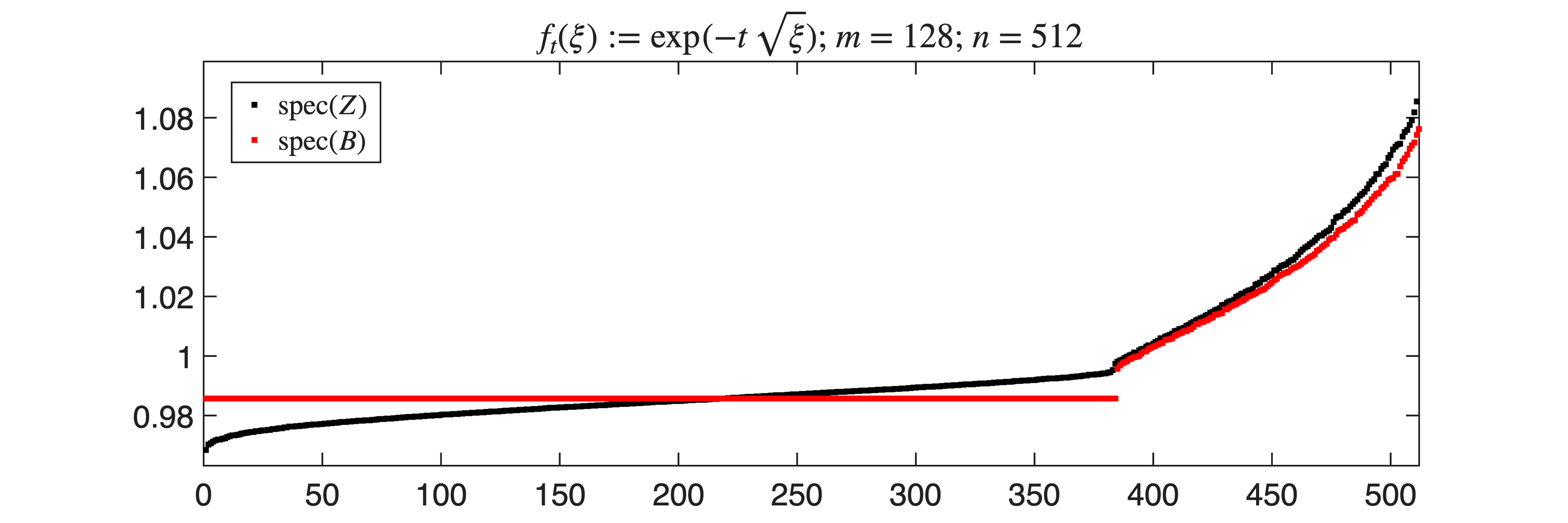}
  \caption{The spectrum of $Z(t_+)$ (black) and {\color{red}the limiting spectral distribution from \eqref{eq:sameSpectralDistributionAsRKM} (red)}. The largest eigenvalue $\approx 3.15$ of $Z(t_+)$ is not shown.}
  \label{fig:random_kernel_matrix_m128n512}
\end{figure}

\begin{figure}[htbp]
  \centering
  \includegraphics[trim = 0mm 0mm 0mm 0mm, clip, width=\textwidth,keepaspectratio]{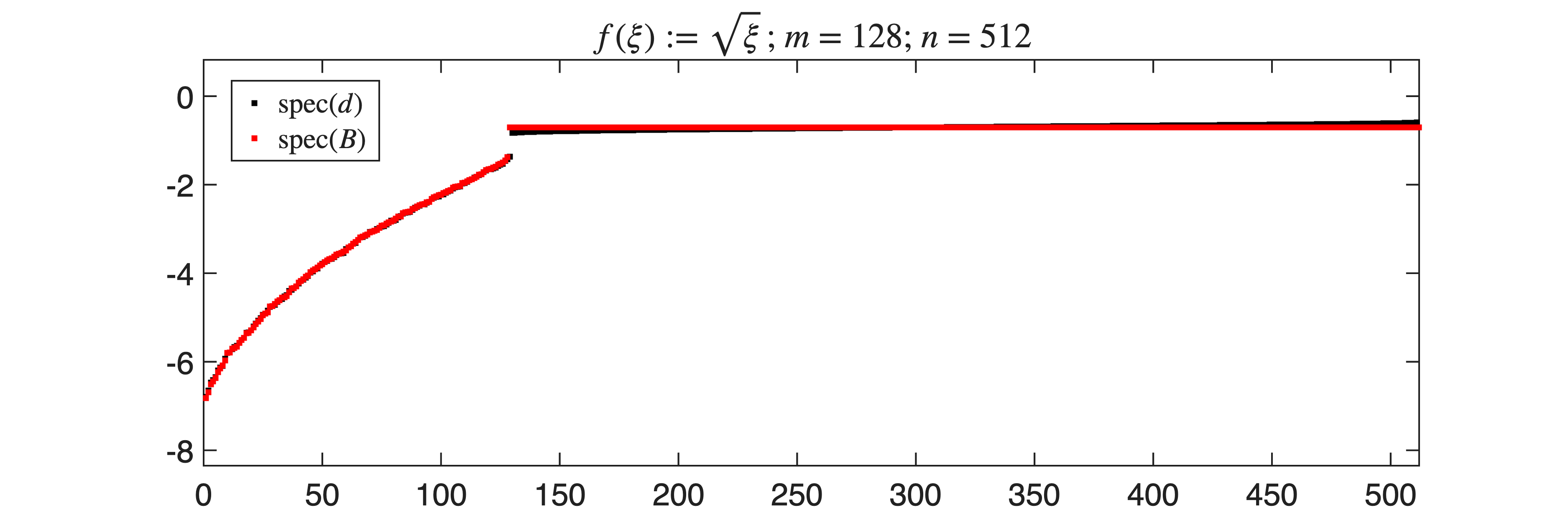}
  \caption{The spectrum of $d$ (black) and {\color{red}the limiting spectral distribution from \eqref{eq:sameSpectralDistributionAsRKM} (red)}. The largest eigenvalue $\approx 720$ of $d$ is not shown.}
  \label{fig:random_kernel_matrix_m128n512_scale0}
\end{figure}

In our case, we have $f_t(\xi) := \exp \left ( -t \sqrt{\xi} \right)$, and example spectra of $Z(t_+)$ and $B$ are shown in Figure \ref{fig:random_kernel_matrix_m128n512}. For $t = 0$, we have instead $f(\xi) = \sqrt{\xi}$, but the same theorem still applies, as shown in Figure \ref{fig:random_kernel_matrix_m128n512_scale0}.

Unfortunately, it is not clear if or how we can usefully leverage knowledge of the limiting spectral distribution in general. For example, there is probably little value in constructing a preconditioner to accelerate computations on normally distributed data. While the relevant theorem also applies to normalized uniform samples from certain nice Riemannian submanifolds with the canonical metric inherited from $\mathbb{R}^m$ (e.g., spheres) \citep{karoui2010spectrum}, the confluence of uniformity and metric requirements render this ostensibly broader applicability largely irrelevant in practice. Furthermore, convergence to the limiting distribution is poor unless $m \gg 1$, in which case the curse of dimensionality suggests that any sample suitable for the theorem would be hard to distinguish than a sample from the standard normal distribution.

\subsection{\label{sec:semianalyticalSummary}Summary}

The semianalytical approximations for scale zero weightings that we have considered all fail. While cutoff scales for isotropic Gaussian-distributed points are easy to approximate using \eqref{eq:gaussianCutoff} and the relation $t_+(\lambda \cdot d) = \lambda^{-1} \cdot t_+(d)$, the utility of this approximation appears to be insignificant in practice. Similarly, it is not clear if how to obtain an estimate of the spectrum of $Z$ or of $d$ outside of situations where its utility is basically zero (e.g., a preconditioner for normally distributed data). In short, semianalytical techniques appear to have little or no practical value for solving the weighting equation, even approximately.

\end{document}